\documentclass{article}

\usepackage{amsfonts}
\usepackage{amssymb}
\usepackage{enumerate}
\usepackage{amsmath}
\usepackage[T1]{fontenc}
\usepackage{graphicx}

\newtheorem{ttt}{Theorem}[section]
\newtheorem{llll}[ttt]{Lemma}
\newtheorem{ccc}[ttt]{Claim}
\newtheorem{eee}[ttt]{Example}
\newtheorem{fff}[ttt]{Fact}
\newtheorem{rrr}[ttt]{Remark}
\newtheorem{sss}[ttt]{Statement}
\newtheorem{ddd}[ttt]{Definition}
\newtheorem{qqq}[ttt]{Question}
\newtheorem{cccc}[ttt]{Corollary}
\newtheorem{nnn}[ttt]{Notation}
\newtheorem{ppp}[ttt]{Problem}
\newtheorem{ccccc}[ttt]{Conjecture}

\newcommand{\beq}{\begin{equation} }

\newcommand{\bt}{\begin{ttt}}
\newcommand{\bl}{\begin{llll}}
\newcommand{\bc}{\begin{ccc}}
\newcommand{\bex}{\begin{eee}}
\newcommand{\bfa}{\begin{fff}}
\newcommand{\br}{\begin{rrr}\upshape}
\newcommand{\bst}{\begin{sss}}
\newcommand{\bd}{\begin{ddd}\upshape}
\newcommand{\bq}{\begin{qqq}}
\newcommand{\bnn}{\begin{nnn}}
\newcommand{\bpr}{\begin{ppp}}
\newcommand{\bcor}{\begin{cccc}}
\newcommand{\bcon}{\begin{ccccc}}

\newcommand{\eeq}{\end{equation}}
\newcommand{\et}{\end{ttt}}
\newcommand{\el}{\end{llll}}
\newcommand{\ec}{\end{ccc}}
\newcommand{\eex}{\end{eee}}
\newcommand{\efa}{\end{fff}}
\newcommand{\er}{\end{rrr}}
\newcommand{\est}{\end{sss}}
\newcommand{\ed}{\end{ddd}}
\newcommand{\eq}{\end{qqq}}
\newcommand{\ecor}{\end{cccc}}
\newcommand{\econ}{\end{ccccc}}
\newcommand{\enn}{\end{nnn}}
\newcommand{\epr}{\end{ppp}}

\newcommand{\bp}{\noindent\textbf{Proof. }}
\newcommand{\ep}{\hspace{\stretch{1}}$\square$\medskip}

\newcommand{\lab}[1]{\label{#1}}

\newcommand{\BB}{\mathbb{B}}
\newcommand{\NN}{\mathbb{N}}

\newcommand{\PP}{\mathbb{P}}

\newcommand{\RR}{\mathbb{R}}

\newcommand{\e}{\varepsilon}

\newcommand{\om}{\omega}
\newcommand{\si}{\sigma}
\newcommand{\ka}{\kappa}
\newcommand{\la}{\lambda}

\newcommand{\iA}{\mathcal{A}}

\newcommand{\iI}{\mathcal{I}}

\newcommand{\iM}{\mathcal{M}}

\newcommand{\iS}{\mathcal{S}}

\newcommand{\iZ}{\mathcal{Z}}

\newcommand{\sm}{\setminus}

\newcommand{\beeq}{\begin{equation}}
\newcommand{\eeeq}{\end{equation}}
\def\su{\subset}

\numberwithin{equation}{section}

\title{A covering theorem and the random-indestructibility of the density zero ideal}

\author{
M\'arton Elekes\thanks{Partially supported by
    Hungarian Scientific Foundation grants no.~83726, 61600, 72655, and J\'anos Bolyai Fellowship.}\\
\\
Alfr\'ed R\'enyi Institute of Mathematics\\
Hungarian Academy of Sciences\\
P.O. Box 127, H-1364 Budapest, Hungary\\
emarci@renyi.hu\\
www.renyi.hu/ $\tilde{}$ emarci\\
and\\
E\"otv\"os Lor\'and University\\
Department of Analysis\\
P\'az\-m\'any P. s. 1/c, H-1117, Budapest, Hungary\\}

\begin{document}

\maketitle

%\subjclass[2000]{Primary 26A99 Secondary 39B22, 39B52, 39B72, 51M99}
\insert\footins{\footnotesize{MSC codes: Primary 03E40 Secondary
    28A99, 03E35}}
\insert\footins{\footnotesize{Key Words: density zero, analytic, P-ideal,
random, forcing, pseudo-intersection, indestructible}}

\begin{abstract}
The main goal of this note is to prove the following theorem. 
If $A_n$ is a sequence of measurable sets in a $\si$-finite
measure space $(X, \iA, \mu)$ that covers $\mu$-a.e. $x \in X$ infinitely many
times, then there exists a sequence of integers $n_i$ of density zero so that $A_{n_i}$
still covers $\mu$-a.e. $x \in X$ infinitely many times. The proof is a probabilistic construction.

As an application we give a simple direct proof of the known theorem that
the ideal of density zero subsets of
the natural numbers is random-indestructible, that is, random forcing does not
add a co-infinite set of naturals that almost contains every ground model
density zero set. This answers a question of B. Farkas.
\end{abstract}

\section{Introduction}

Maximal almost disjoint (MAD) families of subsets of the naturals play
a central role in set theory. (Two sets are \emph{almost disjoint} if
there intersection is finite.) A fundamental question is whether MAD families
remain maximal in forcing extensions. This is often studied in a
little more generality as follows. For a MAD family $\iM$ let
$\iI_\iM$ be the ideal of sets that can be almost contained in a finite
union of members of $\iM$. (\emph{Almost contained} means that only
finitely many elements are not contained.) Then it is easy to see that
$\iM$ remains MAD in a forcing extension if and only if there is no
co-infinite set of naturals in the extension
that almost contains every (ground model) member of $\iI_\iM$. Hence
the following definition is natural.

\bd
An ideal $\iI$ of subsets of the naturals is called \emph{tall} if there is no
co-infinite set that almost contains every member of $\iI$. Let $\iI$
be a tall ideal and $\PP$ be a forcing notion. We say that $\iI$ is
\emph{$\PP$-indestructible} if $\iI$ remains tall after forcing with $\PP$.
\ed

This notion is thoroughly investigated for various well-known ideals
and forcing notions, for instance Hern\'andez-Hern\'andez and Hru\v s\'ak 
proved that the ideal of density zero subsets (see. Definition 
\ref{d:Z}) of the natural numbers is random-indestructible. (Indeed,
just combine \cite[Thm 3.14]{HH}, which is a result of Brendle and
Yatabe, and \cite[Thm 3.4]{HH}.)
B.~Farkas asked if there is a simple and direct proof of this fact.
In this note we provide such a proof. 

This proof actually led us to a covering theorem (Thm. \ref{t:cover}) which we find
very interesting in its own right from the measure theory point of view.
First we prove this theorem in Section \ref{s:cover} by a
probabilistic argument, then we apply it in Section \ref{s:indest} to
reprove that the density zero ideal is random-indestructible
(Corollary \ref{c:indest}), and
finally we pose some problems in Section \ref{s:problems}.

\section{A covering theorem}
\lab{s:cover}

Cardinality of a set $A$ is denoted by $|A|$.

\bd
\lab{d:Z}
A set $A \su \NN$ is of \emph{density zero} if $\lim_{n \to \infty} \frac{ |A
\cap \{0, \dots, n-1 \}  | }{n} = 0$. The ideal of density zero sets is
denoted by $\iZ$. 
\ed

$A \su^* B$ means that $B$ \emph{almost contains} $A$, that is, $A \sm B$ is
finite. The following is well-known.

\bfa
$\iZ$ is a P-ideal, that is, for every sequence $Z_n \in \iZ$ there exists $Z
\in \iZ$ so that $Z_n \su^* Z$ for every $n \in \NN$.
\efa

\bl
\lab{l:everywhere}
Let $(X,\iA, \mu)$ be a measure space of $\si$-finite measure, and let
$\{A_n\}_{n \in \NN}$ be a sequence of measurable sets. Suppose that there
exists $0 = N_0 < N_1 < N_2 < \dots$ so that $A_{N_{k-1}}, \dots, A_{N_k - 1}$
is a cover of $X$ for every $k \in \NN^+$, and also that $k$ divides $N_k -
N_{k-1}$ for every $k \in \NN^+$. Then there exists a set $Z \in \iZ$ 
so that $\{A_n\}_{n \in Z}$ covers $\mu$-a.e. every $x \in X$
infinitely many times.
\el

\bp
Write $\{N_{k-1}, \dots, N_k-1 \} = W^k_0 \cup \cdots \cup W^k_{k-1}$, where the
$W^k_i$'s are the $k$ disjoint arithmetic progressions of difference $k$.
Let $\{\xi_k\}_{k \in \NN^+}$ be a sequence of independent random variables so that
$\xi_k$ is uniformly distributed on $\{0, \dots, k-1 \}$. Define
\[
Z = \cup_{k \in \NN^+} W^k_{\xi_k}.
\]
It is easy to see that $Z \in \iZ$. Hence it suffices to show
that with probability 1 $\mu$-a.e.~$x \in X$ is covered infinitely many
times by $\{A_n\}_{n \in Z}$.

Let us now fix an $x \in X$. Let $E_k$ be the event $\{x \in \cup_{n \in
W^k_{\xi_k}} A_n \}$, that is, $x$ is covered by the set chosen in the $k^{th}$
block. As the $k^{th}$ block is a cover of $X$, $Pr(E_k)
\ge \frac1k$, so $\sum_{k \in \NN^+} Pr(E_k) = \infty$. Moreover, the events $\{
E_k \}_{k \in \NN^+}$ are independent. Hence by the second Borel-Cantelli Lemma
$Pr(\textrm{Infinitely many of the } E_k\textrm{'s occur}) =
1$. So every fixed $x$ is covered infinitely many times with probability 1,
but then by the Fubini theorem with probability 1 $\mu$-a.e. $x$ is covered infinitely
many times, and we are done. (To be more precise, let $(\Omega, \iS, Pr)$ be
the probability measure space, then $Z(\om) = \cup_{k \in \NN}
W^k_{\xi_k(\om)}$. Since the sets $\{ (x,\om) : x \in A_n\}$ and $\{ (x,\om) :
\xi_k(\om) = n\}$ are clearly $\iA \times \iS$-measurable, it is
straightforward to show that
\[
\{ (x, \om) : x \textrm{ is covered infinitely many
  times by } \{ A_n \}_{ n \in Z(\om)} \} \su X \times \Omega
\]
is $\iA \times \iS$-measurable, and hence Fubini applies.)
\ep

\bl
\lab{l:finite}
Let $(X,\iA, \mu)$ be a measure space of finite measure, and let $\{A_n\}_{n
\in \NN}$ be a sequence of measurable sets that covers $\mu$-a.e. every $x
\in X$ infinitely many times. Then there exists a set $Z \in \iZ$ so
that $\{A_n\}_{n \in Z}$ still covers $\mu$-a.e. every $x \in X$
infinitely many times.
\el

\bp
Let $\e > 0$ be arbitrary and set $N_0 = 0$. 
By the continuity of measures, there exists $N_1$ so that $\mu
( X \sm (A_{N_0} \cup \dots \cup A_{N_1 -1})) \le \frac\e2$. 
Since $\{A_n\}_{n \ge N_1}$
still covers $\mu$-a.e. $x \in X$ infinitely many times, we can continue this
procedure, and recursively define $0 = N_0 < N_1 < N_2 < \dots$ so that $\mu
(X \sm (A_{N_{k-1}} \cup \dots \cup A_{N_k - 1})) \le \frac\e{2^k}$ for every $k \in
\NN^+$. We can also assume (by choosing larger $N_k$'s at each step) that $k$
divides $N_k - N_{k-1}$ for every $k \in \NN^+$.

Let $X_\e = \cap_{k \in \NN^+} (A_{N_{k-1}} \cup \dots \cup A_{N_k - 1})$, then
$\mu (X \sm X_\e) \le \e$. Let us restrict $\iA$, the $A_n$'s and $\mu$
to $X_\e$, and apply the previous lemma with this setup to obtain $Z_\e$.

Let us now consider $\e = 1, \frac12, \frac13, \dots$, then for every $m \in
\NN^+$ every $x \in X_\frac1m$ is covered infinitely many times by $\{ A_n
\}_{n \in Z_\frac1m}$. Since $\iZ$ is a P-ideal, there exists a $Z \in \iZ$
such that $Z_\frac1m \su^* Z$ for every $m$. Hence for every $m \in
\NN^+$ every $x \in X_\frac1m$ is covered infinitely many times by $\{ A_n
\}_{n \in Z}$. But then we are done, since $\mu$-a.e. $x
\in X$ is in $\cup_m X_\frac1m$. 
\ep

\bt
\lab{t:cover}
Let $(X,\iA, \mu)$ be a measure space of $\si$-finite measure, and let
$\{A_n\}_{n \in \NN}$ be a sequence of measurable sets that covers
$\mu$-a.e. every $x \in X$ infinitely many times. Then there exists a set $Z
\su \NN$ of density zero so that $\{A_n\}_{n \in Z}$ still covers
$\mu$-a.e. every $x \in X$ infinitely many times.
\et

\bp
Write $X = \cup X_m$, where each $X_m$ is of finite measure. For each $X_m$
obtain $Z_m$ by the previous lemma. Then a $Z \in \iZ$ such that $Z_m \su^* Z$
for every $m$ clearly works.
\ep

The following example shows that the purely topological analogue of Theorem
\ref{t:cover} is false.

\bex
There exists a sequence $U_n$ of clopen sets covering every point of
the Cantor space infinitely many times so that for
every $Z \in \iZ$ there exists a point covered only finitely many times by $\{
U_n : n \in Z \}$.
\eex

\bp
By an easy recursion we can define a sequence $U_n$ of clopen subsets of the
Cantor set $C$ and a sequence of naturals $0=N_0<N_1<\dots$ with the following
properties.
\begin{enumerate}
\item
$U_{N_{k-1}}, \dots, U_{N_k-1}$ (called a `block') is a disjoint cover of $C$,
\item
every block is a refinement of the previous one,
\item
\lab{3}
if $U_n$ is in the $k^{th}$ block and is partitioned into $U_t, \dots , U_s$
in the $k+1^{st}$ block (called the `immediate successors of $U_n$') then $s
\ge 2t$.
\end{enumerate}
Let $Z \in \iZ$ be given, and let $n_0$ be so that $\frac{|Z \cap \{0, \dots,
  n-1\}|}{n} < \frac12$ for every $n \ge n_0$. By \ref{3}.~$\{ U_n : n \in
Z \}$ cannot contain all immediate successors of any $U_m$ above
$n_0$. Therefore, starting at a far enough block, we can recursively pick a
$U_{n_i}$ from each block so that
$n_i \notin Z$ for every $i$, and $\{ U_{n_i} \}_{i \in \NN}$ 
is a nested sequence of clopen
sets. But then the intersection of this sequence is only covered finitely many
times by $\{U_n : n \in Z\}$.
\ep

\br
We can `embed' this example into any topological space containing a copy
of the Cantor set (e.g. to any uncountable Polish space) by just adding
the complement of the Cantor set to all
$U_n$'s. Of course, the new $U_n$'s will only be open, not clopen.
\er

\section{An application: The density zero ideal is random-indestructible}
\lab{s:indest}

In this section we give a simple and direct proof of the
random-indestructibility of $\iZ$, which was first proved in
\cite{HH}.

$[\NN]^\om$ denotes the set of infinite subsets of $\NN$. Since it can be 
identified with a $G_\delta$ subspace of $2^\omega$ in the natural way,
it carries a Polish space topology where the sub-basic open sets are the
sets of the form $[n] = \{A \in [\NN]^\om : n \in A\}$ and their
complements. Let $\la$ denote Lebesgue measure.

\bl
For every Borel function $f : \RR \to [\NN]^\om$ there exists a set $Z \in
\iZ$ such that $f(x) \cap Z$ is infinite for $\la$-a.e. $x \in \RR$.   
\el

\bp
Let $A_n = f^{-1}([n])$, then $A_n$ is clearly Borel, hence Lebesgue
measurable. For every $x \in \RR$ 
\beq
\lab{e:1}
x \in A_n \iff x \in f^{-1}([n]) \iff f(x) \in [n] \iff n \in f(x).
\eeq 
Since every $f(x)$ is infinite, \eqref{e:1} yields that every $x \in \RR$ is
covered by infinitely many $A_n$'s. By Theorem \ref{t:cover} there exists a $Z
\in \iZ$ such that for $\la$-a.e. $x \in \RR$ we have $x \in A_n$ for infinitely
many $n \in Z$. But then by \eqref{e:1} for $\la$-a.e. $x \in \RR$ we have $n
\in f(x)$ for infinitely many $n \in Z$, so $f(x) \cap Z$ is infinite.
\ep

Recall that \emph{random forcing} is $\BB = \{p \su \RR : p \textrm{ is Borel},
\la(p) > 0 \}$ ordered by inclusion. The \emph{random real $r$} is
defined by $\{r\} = \cap_{p \in
G} \ p$, where $G$ is the generic filter. For the terminology and basic facts
concerning random forcing consult e.g. \cite{Ku}, \cite{Je}, \cite{BJ}, or
\cite{Za}. In particular, we will assume familiarity with coding of Borel sets
and functions, and will freely use the same symbol for all versions of a Borel
set or function.  The following fact is well-known and easy to prove.

\bfa
\lab{f:random}
Let $B \su \RR$ be Borel. Then $p \Vdash ``r \in B"$ iff $\la(p \sm B) = 0$.
\efa

\bcor
\lab{c:indest}
The ideal of density zero subsets of
the natural numbers is random-indestructible, that is, random forcing does not
add a co-infinite set of naturals that almost contains every ground model
density zero set.
\ecor

\bp
For a Borel function $f : \RR \to [\NN]^\om$ and a set $Z \in
\iZ$ let 
\[
B_{f,Z} = \{x \in \RR : f(x) \cap Z \textrm{ is infinite}\},
\]
then by the previous lemma for every $f$ there is a $Z$ so that $B_{f,Z}$ is
of full measure. By Fact \ref{f:random} for every $f$ there is a $Z$ so that
$1_\BB \Vdash ``f(r) \cap Z \textrm{ is infinite}"$. Hence for every $f$
$1_\BB \Vdash ``\exists Z \in \iZ \cap V \textrm{ so that } f(r) \cap Z
\textrm{ is infinite}"$. But every $y \in [\NN]^\om \cap V[r]$ is of the form
$f(r)$ for some ground model (coded) Borel function $f : \RR \to [\NN]^\om$,
so we obtain that for every $y \in [\NN]^\om \cap V[r]$
$1_\BB \Vdash ``\exists Z \in \iZ \cap V \textrm{ so that } y \cap Z
\textrm{ is infinite}"$. Therefore $1_\BB \Vdash ``\forall y \in [\NN]^\om
\exists Z \in \iZ \cap V \textrm{ so that } y \cap Z \textrm{ is infinite}"$,
and setting $x = \NN \sm y$ yields $1_\BB \Vdash ``\forall x \su \om \textrm{
co-infinite } \exists Z \in \iZ \cap V \textrm{ so that } Z \not\su^* x$'', so
we are done.
\ep

\br
Clearly, $\iZ$ is also $\BB(\kappa)$-indestructible, since every new real is
already added by sub-poset isomorphic to $\BB$. ($\BB(\kappa)$ is the
usual poset for adding $\kappa$ many random reals by the measure algebra on
$2^\ka$.)
\er

\br
\lab{r:DC}
The referee of this paper has pointed out that all arguments of the paper can 
actually be carried out in the axiom system $ZF + DC$. ($ZF$ is the usual 
Zermelo-Fraenkel axiom system without the Axiom of Choice, and $DC$ is the Axiom
of Dependent Choice.) Hence Corollary \ref{c:indest} actually applies to
forcing over a model of $ZF + DC$ as well.
\er

\section{Problems}
\lab{s:problems}

There are numerous natural directions in which one can ask questions in light
of Corollary \ref{c:indest} and Theorem \ref{t:cover}. As for the former one,
one can consult e.g. \cite{BY} and the references therein. As for the latter one,
it would be interesting to investigate what happens if we replace the density
zero ideal by another well-known one, or if we
replace the measure setup by the Baire category analogue, or if we consider
non-negative functions (summing up to infinity a.e.) instead of sets, or even if we
consider $\kappa$-fold covers and ideals on $\kappa$.

\bigskip
\noindent
\textbf{Acknowledgment} The author is indebted to Rich\'ard Balka for some 
helpful discussions and to an anonymous referee for suggesting Remark \ref{r:DC}.


\begin{thebibliography}{1}

%\bibitem{BE} R.~Balka, M.~Elekes, The structure of rigid functions,
%\textsl{J.~Math.~Anal.~Appl.} \textbf{345}, no.~2, (2008), 880--888.

\bibitem{BJ} T.~Bartoszy\'nski, H.~Judah, \textsl{Set theory. On the structure
of the real line.} A K Peters, Ltd., Wellesley, MA, 1995.

\bibitem{BY} J.~Brendle, S.~Yatabe, Forcing indestructibility of MAD families,
\textsl{Ann. Pure Appl. Logic} \textbf{132} (2005), no. 2-3, 271--312. 

%\bibitem{CCR} B.~Cain, J.~Clark, D.~Rose, Vertically rigid functions,
%\textsl{Real Anal.~Exchange} \textbf{31}, no.~2, (2005/2006), 515--518.

%\bibitem{Fa} K.~J.~Falconer, \textsl{The geometry of fractal sets.}
%Cambridge Tracts in Mathematics No. 85, Cambridge University
%Press, 1986.

\bibitem{HH} F.~Hern\'andez-Hern\'andez, M.~Hru\v s\'ak,
Cardinal invariants of analytic $P$-ideals.
\textsl{Can. J. Math.} \textbf{59} (2007), No. 3, 575-595. 

%\bibitem{HG} M.~Hru\v s\'ak, S.~Garc\'\i a Ferreira,
%Ordering MAD families a la Katétov. \textsl{J. Symbolic Logic} \textbf{68}
%(2003), no.~4, 1337--1353.

\bibitem{Je} T.~Jech, \textsl{Set theory. The third millennium edition,
revised and expanded.} Springer Monographs in Mathematics. Springer-Verlag,
Berlin, 2003.

%\bibitem{Ke} A.~S.~Kechris, \textsl{Classical Descriptive Set Theory.}
%Springer-Verlag, 1995.

%\bibitem{Ku} K.~Kuratowski: \textsl{Topology.} Academic Press, 1966.

\bibitem{Ku} Kunen, K.: \textsl{Set theory. An introduction to independence
proofs.} Studies in Logic and the Foundations of Mathematics,
102. North-Holland, 1980.

%\bibitem{Ma} P.~Mattila: \textsl{Geometry of Sets and Measures in
%Euclidean Spaces.} Cambridge Studies in Advanced Mathematics
%No.~44, Cambridge University Press, 1995.

%\bibitem{Ox} J.~C.~Oxtoby: \textsl{Measure and Category. A survey of the
%analogies between topological and measure spaces.} Second edition.
%Graduate Texts in Mathematics No. 2, Springer-Verlag, 1980.

%\bibitem{Ri} C.~Richter, Continuous rigid functions, preprint, see
%http://www.minet.uni-jena.de/Math-Net/reports/shadows//08-02report.html

%\bibitem{vN} J.~von Neumann, Ein System algebraisch unabh\"angiger Zahlen,
%\textsl{Math. Ann.} \textbf{99}, (1928), 134-141.

\bibitem{Za} J.~Zapletal, \textsl{Forcing idealized.} Cambridge Tracts in
Mathematics, 174. Cambridge University Press, Cambridge, 2008. 

\end{thebibliography}
\end{document}